\documentclass{amsart}
\usepackage{amsfonts}
\usepackage{amsmath,amssymb}
\usepackage{amsthm}
\usepackage{amscd}
\usepackage{graphics}
\usepackage{graphicx}

\theoremstyle{remark}{

\newtheorem{Prob}{{\rm Problem}}

}
\theoremstyle{plain}
{

\newtheorem{Thm}{Theorem}

}

\begin{document}
\title[A new classification of Morse functions on $3$-dimensional manifolds]{On a classification of Morse functions on $3$-dimensional manifolds represented as connected sums of manifolds of Heegaard genus one}
\author{Naoki kitazawa}

\address{Institute of Mathematics for Industry, Kyushu University, 744 Motooka, Nishi-ku Fukuoka 819-0395, Japan\\
 TEL (Office): +81-92-802-4402 \\
 FAX (Office): +81-92-802-4405 \\
}
\email{n-kitazawa@imi.kyushu-u.ac.jp, naokikitazawa.formath@gmail.com}
\urladdr{https://naokikitazawa.github.io/NaokiKitazawa.html}
\maketitle
\begin{abstract}
Morse functions are important objects and tools in understanding topologies of manifolds since the 20th century. Their classification has been natural and difficult problems, and surprisingly, this is recently developing. Since the 2010's, results for cases of surfaces have been presented by Gelbukh, Marzantowicz and Michalak for example. We have also longed for higher dimensional cases. We present a classification of Morse functions on $3$-dimensional manifolds represented as connected sums of manifolds of Heegaard genus one. We concentrate on Morse functions such that preimages of single points containing no singular points are disjoint unions of spheres and tori. Existence of such functions implies that the $3$-dimensional closed and connected manifolds are of such manifolds. It has been shown by Saeki in 2006 and we further study structures of these functions.

\end{abstract}
\section{Introduction.}
\label{sec:1}
Morse functions have been fundamental and important tools in differential topology of manifolds. They themselves are also interesting objects.
Related theory dates back to the former half of the 20th century. 

As an explicit study, classifications of Morse functions on given manifolds (of certain classes) have been natural and important studies. Surprisingly, such studies are very new and developing recently. Since the 2010s, related studies have been presented by \cite{gelbukh1, gelbukh2, marzantowiczmichalak, michalak1, michalak2} for example. Most of them are for closed surfaces. These studies also long for higher dimensional cases. We present several important and motivating studies. Before this, we review several notions for example. 

Let ${\mathbb{R}}^k$ denote the $k$-dimensional Euclidean space with $\mathbb{R}={\mathbb{R}}^1$.
For a smooth map $c:X \rightarrow Y$ between smooth manifolds $X$ and $Y$, a point $p \in X$ is a singular point of $c$ if the rank of the differential drops there. A Morse function $c:X \rightarrow \mathbb{R}$ means a smooth 
function such that each singular point is in the interior of the manifold $X$ and that each singular point $p$ is of the form $c(x_1,\cdots x_m)={\Sigma}_{j=1}^{m-i(p)} {x_j}^2-{\Sigma}_{j=1}^{i(p)} {x_{m-i(p)+1}}^2+c(p)$ for some integer $0 \leq i(p) \leq m$. The integer $i(p)$ is shown to be uniquely defined and the {\it index} of $p$ {\it for $c$}. See \cite{milnor} for theory of Morse functions and applications to manifold theory, for example. Here, graphs are fundamental and important. They are CW complexes consisting of 1-cells ({\it edges}) and 0-cells ({\it vertices}). 
The {\it edge set} of the graph is the set of all edges. The {\it vertex set} of the graph is the set of all vertices. We can orient edges: the graph is a {\it digraph}. Two graphs are {\it isomorphic} if there exists a piecewise smooth homeomorphism between them mapping the vertex set of a graph onto that of another graph: this map is an {\it isomorphism} between the graphs.
For two digraphs, they are {\it isomorphic} if they are isomorphic as graphs, admitting an isomorphism of the graphs mapping each edge into another edge by preserving the orientations. Hereafter, the {\it degree} of a vertex is the number of edges containing the vertex.
 
The {\it Reeb space} $W_c$ of a smooth function $c:X \rightarrow \mathbb{R}$ is defined as the quotient space $X/{\sim}_c$ by the equivalence relation; $x_1 {\sim}_c x_2$ if and only if $x_1$ and $x_2$ are in a connected component of a same preimage $c^{-1}(y)$. We can define the quotient map $q_c:X \rightarrow W_c$ and in a unique way we can define a function $\bar{c}:W_c \rightarrow \mathbb{R}$ by the relation $c=\bar{c} \circ q_c$. By defining the vertex set of $W_c$ as the set of all connected components containing some singular points of $c:X \rightarrow \mathbb{R}$ (where $X$ is with no boundary), this is a graph in considerable cases and the {\it Reeb graph} of $c$: more rigorously, \cite{saeki2} discusses this rigorously and our cases are of such cases unless otherwise stated.
We can also see the Reeb graph $W_c$ as a digraph in a canonical way. We orient an edge connecting $v_1$ and $v_2$ as an edge starting from $v_1$ to $v_2$ if $\bar{c}(v_1)<\bar{c}(v_2)$. Note that the function $\bar{c}$ is piecewise smooth and injective at each edge from the definition. We also call $W_c$ the {\it Reeb digraph} of $c$. We can orient a graph $K$ in the same way if there exists a piecewise smooth function $g:K \rightarrow \mathbb{R}$ such that at each edge it is injective. Let $K_g$ denote the digraph.

We go back to some theorems. We omit precise exposition on fundamental notions such as the {\it genus} of a closed and connected surface, properties, and theorems on closed surfaces. 

\begin{Thm}[\cite{michalak1}]
\label{thm:1}
Let a finite and connected digraph $K_g$ be defined from a piecewise smooth function $g:K \rightarrow \mathbb{R}$ with the following conditions.  
\begin{enumerate}
\item The graph $K_g$ is not isomorphic to a connected graph with exactly one edge and two vertices{\rm :} we do not consider Morse functions with exactly two singular points on spheres for so-called Reeb's theorem {\rm (}\cite{reeb}{\rm )}.
\item The 1st Betti number of the graph $K$ is $a$. The number of vertices of degree $2$ is $b$.
\item If $g$ has a local extremum at a vertex $v \in K$, then $v$ is of degree $1$.
\end{enumerate} 
Then for a closed, orientable and connected surface $M$ of genus $g \geq a+b$ or a closed, connected and non-orientable surface $M$ of genus $2a+b$, there exists a Morse function $f:M \rightarrow \mathbb{R}$ whose Reeb digraph $W_f$ and the digraph $K_g$ are isomorphic. 
Conversely, for such a Reeb graph, if the domain of the Morse function is a closed and connected surface, then the surface is diffeomorphic to one of these surfaces. 
\end{Thm}

\cite{gelbukh4} is an extended result for Morse-Bott functions on closed surfaces. A {\it Morse-Bott} function is a generalized function where at each singular point the function is represented as the composition of a projection (, a smooth map with no singular points,) with a Morse function. See the third section. 

Furthermore, \cite{marzantowiczmichalak} studies systems of smooth hypersurfaces in compact and connected manifolds and homomorphisms of the fundamental groups onto the free groups induced by the systems naturally. This also says that such homomorphisms are induced by Morse functions and their Reeb graphs. This also studies deformations and some classifications of Morse functions on closed surfaces and higher dimensional manifolds with \cite{michalak2}. It also implies that precise classifications for higher dimensional cases respecting preimages are hard. \cite{michalak3} studies fundamental properties on relations between fundamental groups of $3$-dimensional closed and connected manifolds and numbers of singular points of Morse functions on them. Explicit cases for circle bundles over closed, orientable and connected surfaces are also studied. The study of Michalak is independent of and closely related to our study.

Our main result is Theorem \ref{thm:2}, a kind of higher dimensional versions of Theorem \ref{thm:1}. See \cite{hempel} for example for $3$-dimensional manifold theory. The {\it Heegaard genus} of a $3$-dimensional manifold is a fundamental and important numerical invariant for $3$-dimensional closed and connected manifolds. Explicitly, a manifold diffeomorphic to the $3$-dimensional sphere $S^3$ is only the manifold of Heegaard genus $0$. The product $S^2 \times S^1$ of the $2$-dimensional sphere $S^2$ and the circle $S^1$ or a so-called {\it Lens space} is of Heegaard genus $1$ and manifolds of Heegaard genus $1$ must be diffeomorphic to them. Circle bundles over the sphere $S^2$ (which are not trivial bundles or $S^3$, seen as a so-called Hopf-fibration,) are of Lens spaces. 
\begin{Thm}[Our main result]
\label{thm:2}
Let a finite and connected digraph $K_g$ be defined from a piecewise smooth function $g:K \rightarrow \mathbb{R}$ with the following conditions.  
\begin{enumerate}
\item The graph $K_g$ is not isomorphic to a connected graph with exactly one edge and two vertices.
\item The 1st Betti number of the graph $K$ is $a$.
\item If $g$ has a local extremum at a vertex $v \in K$, then $v$ is of degree $1$.
\item A map $l_K:E_K \rightarrow \{0,1\}$ on the edge set $E_K$ of $K$ satisfies that the number of vertices $v$ of degree $2$ contained in two edges $e_{v,1}$ and $e_{v,2}$ with $l_K(e_{v,1})=l_K(e_{v,2})=0$ is $b$ and that the size of the set ${l_K}^{-1}(1)$ is $c$.

\end{enumerate}
Then on a $3$-dimensional closed manifold $M$ diffeomorphic to a manifold represented as a connected sum of $r \geq a+b$ copies of $S^1 \times S^2$ and $c^{\prime} \leq c$ manifolds whose free groups are finite and whose Heegaard genera are $1$, there exists a Morse function $f:M \rightarrow \mathbb{R}$ enjoying the following properties. 
\begin{enumerate}
\item There exists an isomorphism ${\phi}_{f,g}:K_g \rightarrow W_f$ between these digraphs.
\item At an edge $e \in E_K$ with $l_K(e)=0$ and at each point $p_e$ in the interior of the edge ${\phi}_{f,g}(e)$, the preimage ${q_f}^{-1}(p_e)$ is diffeomorphic to $S^2$. At an edge $e \in E_K$ with $l_K(e)=1$ and at each point $p_e$ in the interior of the edge ${\phi}_{f,g}(e)$, the preimage ${q_f}^{-1}(p_e)$ is diffeomorphic to the torus $S^1 \times S^1$.
\end{enumerate}
Conversely, for a Morse function $f:M \rightarrow \mathbb{R}$ on a $3$-dimensional closed and orientable manifold $M$ with such a Reeb digraph $W_f$, $M$ is diffeomorphic to one of these manifolds. 
Here the case $r=0$ and $c^{\prime}=0$ implies $M$ is diffeomorphic to $S^3$. 
\end{Thm}
Together with the presented studies, the following is another motivating study.
\begin{Thm}[{\cite[Theorem 6.5]{saeki1}}]
\label{thm:3}
If a $3$-dimensional closed, connected and orientable manifold admits a Morse function such that preimages of single points containing no singular points of it  are disjoint unions of copies of $S^2$ and $S^1 \times S^1$, then the manifold is as in Theorem \ref{thm:2}. Such manifolds also admit such Morse functions. 
\end{Thm}

For related studies, we refer to reconstruction of smooth functions whose Reeb digraphs and given digraphs are isomorphic. This is started in \cite{sharko}. \cite{gelbukh1, gelbukh2, martinezalfaromezasarmientooliveira, masumotosaeki, michalak1} are of important studies. The author has considered not only graphs, but also information on preimages beforehand in \cite{kitazawa3, kitazawa4, kitazawa5}, as a pioneer.   

We present our proof of Theorem \ref{thm:2} and related arguments in the next section. The third section is for several remarks.
\section{Our main result: Theorem \ref{thm:2}.}

We prove Theorem \ref{thm:2}. We apply arguments from articles presented in the previous section. 
Hereafter, let $D^k:=\{x \in {\mathbb{R}}^k \mid {\Sigma}_{j=1}^{k} {x_j}^2 \leq 1\}$ be the $k$-dimensional unit disk. For a topological manifold $X$, let $\partial X$ denote the boundary of $X$. For example, the ($k-1$)-dimensional unit sphere $S^{k-1}$ is $\partial D^k$.

As a fundamental theory, we also apply a natural one-to-one correspondence between a singular point of index $k$ for a Morse function on an $m$-dimensional manifold and a so-called {\rm (}{\it k-}{\rm )}{\it handle}, a cornered smooth manifold diffeomorphic to $D^k \times D^{m-k}$. 

We present a more precise argument on this. First, singular points of Morse functions appear discretely. 
Finitely many handles are attached to boundaries of submanifolds of codimension $0$ in the manifolds of the domains, disjointly.

Hereafter, for two real numbers $t_1<t_2$, let $(t_1,t_2):=\{t \in \mathbb{R} \mid t_1 < t < t_2\}$,$(t_1,t_2]:=\{t \in \mathbb{R} \mid t_1 < t \leq t_2\}$, $[t_1,t_2):=\{t \in \mathbb{R} \mid t_1 \leq t < t_2\}$, and $[t_1,t_2]:=\{t \in \mathbb{R} \mid t_1 \leq t \leq t_2\}$.

We consider a Morse function $c:X \rightarrow \mathbb{R}$ on an $m$-dimensional smooth compact manifold $X$ such that the image is a closed interval $[t_1,t_2]$. Let $t_1<q<t_2$ be a value realized as a value at a singular point of $c$.
Suppose also that $c^{-1}(t_1) \sqcup c^{-1}(t_2)$ and the boundary $\partial X$ coincide and that each preimage $c^{-1}(t_j)$ ($j=1,2$) is a union of connected components of $\partial X$. Let $\epsilon>0$ be sufficiently small and $c^{-1}([t_1,q+\epsilon))$ is a smooth compact manifold obtained by attaching finitely many handles to $c^{-1}(q-\epsilon) \subset c^{-1}([t_1,q-\epsilon))$. We attach a $k$-handle $D^k \times D^{m-k}$ by a diffeomorphism from $\partial D^k \times D^{m-k}$ onto a submanifold of $\partial X$ diffeomorphic to $\partial D^k \times D^{m-k}$. More rigorously, they are ($m-1$)-dimensional smooth compact submanifolds and seen as trivial linear bundles over a sphere $\partial D^k=S^{k-1}$ and the handles are attached by isomorphisms of bundles. Each handle corresponds to each singular point of the function. 

Let $d$ be the number of handles and singular points of the function here.
We choose an order on the set of the $d$ singular points of the function. We can also deform the function into another Morse function by a suitable homotopy without changing its singular points and the values on $c^{-1}(\mathbb{R}-(q-\epsilon,q+\epsilon))$ in such a way that at distinct singular points of the resulting new function $c^{\prime}:X \rightarrow \mathbb{R}$ the values are distinct. More precisely, we can change the function in such a way that at the $j$-th singular point of the $d$ singular points ($1 \leq j \leq d$) of the function, the value is $q_j$ with $q-\epsilon<q_j<q+\epsilon$ and that the relation $q_{j_1}<q_{j_2}$ is satisfied for $1 \leq j_1<j_2 \leq d$. We can also change the order on the set $\{s_1,\cdots,s_j\}$ in arbitrary ways (there exist $d!$ cases). For each singular point of the new function, the handle same as that in the original function corresponds to. In other words, we can attach handles one after another to the (original) mutually disjoint submanifolds diffeomorphic to $\partial D^k \times D^{m-k-1}$ of the preimage of $q-\epsilon$, according to the values, increasing from $q-\epsilon$ to $q+\epsilon$.

We can consider a converse argument. For such a Morse function $c^{\prime}:X \rightarrow \mathbb{R}$ such that at distinct singular points of the function the values are distinct, handles are attached to mutually disjoint submanifolds diffeomorphic to $\partial D^k \times D^{m-k-1}$ of the preimage ${c^{\prime}}^{-1}(q_j-\epsilon)$ of $q_j-\epsilon$ one after another. Here we abuse the same notation. We can deform the function into another Morse function $c:X \rightarrow \mathbb{R}$ by a suitable homotopy without changing its singular points and the values on ${c^{\prime}}^{-1}(\mathbb{R}-(q-\epsilon,q+\epsilon))$ in such a way that at all singular points of $c$, the values are same and $q$.

Other than this, we need several elementary and fundamental arguments on Morse functions. We do not present them precisely.

For this theory, see \cite{milnor} for example.

\begin{proof}
STEP 1 On a Morse function $c:X \rightarrow \mathbb{R}$ such that at distinct singular points of $c$ the values are distinct and that preimages of single points containing no singular points of it are disjoint unions of copies of the sphere $S^2$ and the torus $S^1 \times S^1$: local information on the Reeb space (Reeb digraph) and preimages. \\
\ \\
Morse functions of this class are essential in our arguments. FIGURE \ref{fig:1} depicts local information on the Reeb (di)graph and preimages. A blue (red) edge shows an edge at a point in the interior of which the preimage is diffeomorphic to the sphere $S^2$ (resp. the torus $S^1 \times S^1$). A black dot is for a singular point of index $0$ for the Morse function. A green dot is for a singular point of index $1$ for the Morse function. We can obtain the other cases by considering the function $-c:X \rightarrow \mathbb{R}$: the index of a singular point for this new Morse function is turned into $2$ from $1$ or $3$ from $0$. This is fundamental and important in our proof.
\begin{figure}
		\includegraphics[width=80mm,height=27.5mm]{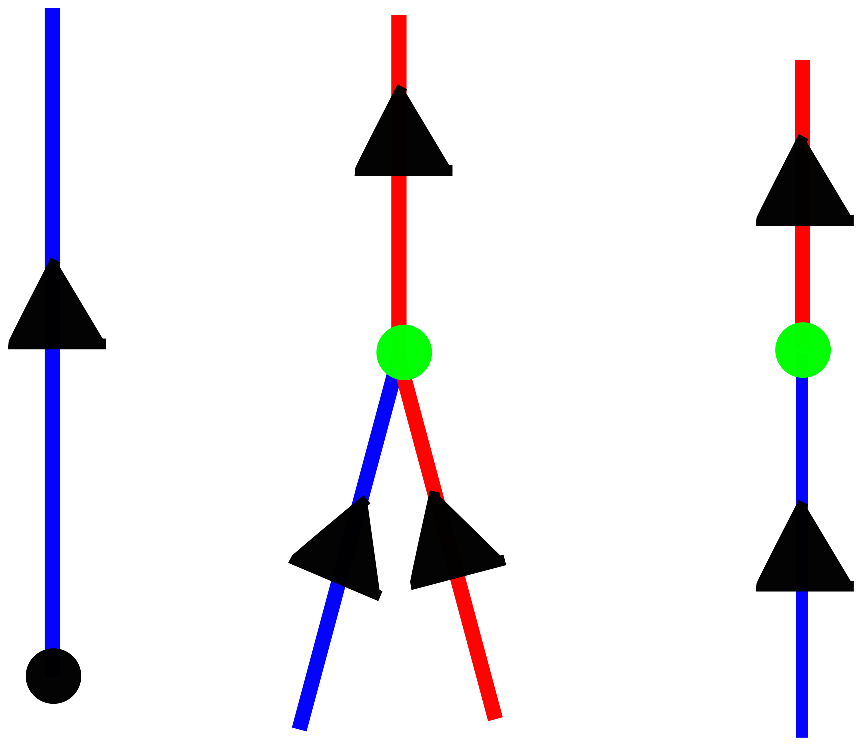}
		\caption{Local information on the Reeb (di)graph and preimages for SSTF Morse functions. A blue (red) edge shows an edge at a point in the interior of which the preimage is diffeomorphic to the sphere $S^2$ (resp. the torus $S^1 \times S^1$). A black (green) dot is for a singular point of index $0$ (resp. 1) for the Morse function.}
\label{fig:1}
	\end{figure}
We can know this as a kind of elementary exercise on Morse functions, singular points of them, and handles.

Hereafter we call such a Morse function a {\it simple sphere-torus-fibered} Morse function or an {\it SSTF} Morse function.

Hereafter we call a Morse function such that at distinct singular points of the function the values are always distinct a {\it simple} Morse function. An SSTF Morse function is simple of course.

Note that the notion of Reeb (di)graph is defined for a smooth function on a manifold with no boundary. However, let us abuse the notion for smooth functions on manifolds with non-empty boundaries
 (where it is no problem). 

\ \\
 STEP 2 (Re)construction of a Morse function $f:M \rightarrow \mathbb{R}$ on a desired $3$-dimensional closed manifold $M$ from a digraph $K_g$: we also obtain an SSTF Morse function $f_0:M \rightarrow \mathbb{R}$ by a suitable homotopy from $f$. \\

This reconstruction is first presented in \cite{michalak1} in the case preimages of single points containing no singular points of the functions are disjoint unions of spheres. Later, in \cite{kitazawa1}, the case preimages of single points containing no singular points of the functions are disjoint unions of (2-dimensional) spheres and closed, connected and orientable surfaces, is studied. 
As more generalized cases, we also respect our studies \cite{kitazawa3, kitazawa4, kitazawa5} for Morse functions on manifolds whose dimensions are greater than $2$ and preimages of single points diffeomorphic to more general manifolds: however we do not assume related knowledge and arguments at all. 

We locally construct a local Morse function around each edge so that the quotient map onto its Reeb space is a map onto a small regular neighborhood of a vertex. After that, we consider a local trivial smooth bundle over each edge whose fiber is as prescribed by the map $l_K$. We glue these local functions together to have our desired Morse function on a manifold.
\begin{enumerate}
\item \label{1}
 First, at each vertex $v$ where $g$ has a local extremum, we consider a function $f_v:X_v \rightarrow \mathbb{R}$ on the $m$-dimensional unit disk $X_v:=D^3 \subset {\mathbb{R}}^m$, represented by $f_v(x)=\pm {\Sigma}_{j=1}^m {x_j}^2 +g(v)$. This is also regarded as an SSTF Morse function $f_{0,v}:X_v \rightarrow \mathbb{R}$.
\item \label{2}
At each vertex $v$ of degree $2$ contained in two distinct edges $e_{v,1}$ and $e_{v,2}$ with $l_K(e_{v,1})=l_K(e_{v,2})=0$, we can locally have a desired Morse function $f_v:X_v \rightarrow \mathbb{R}$ on a $3$-dimensional compact and orientable manifold $X_v$ as follows.
\begin{itemize}
\item At all singular points of the function the values are $g(v)$.
\item The image is a closed interval $[t_{v,1},t_{v,2}]$ containing $g(v)$ in the interior.
\item The preimage of each extremum is diffeomorphic to $S^2$.
\end{itemize}
Conversely, desired Morse functions must be such functions. We can deform such a Morse function by a suitable homotopy to an SSTF Morse function $f_{0,v}:X_v \rightarrow \mathbb{R}$ such that the Reeb space is homeomorphic to a graph of 1st Betti number at least $1$.
 
We explain this more precisely. For this, see also \cite{gelbukh1, gelbukh3, gelbukh4, michalak1} for example.



For a Morse function explained just before, this can be corresponded to some attachment of handles to the sphere
${f_{0,v}}^{-1}(t_{v,1}+\epsilon)$ to obtain $X_v={f_{0,v}}^{-1}([t_{v,1},t_{v,2}])$ according to the following procedure where $\epsilon>0$ is a sufficiently small number.
This respects \cite[Lemma 6.6]{saeki1}, presented as Theorem \ref{thm:5} later, and its original proof. More precisely, in Theorem \ref{thm:5}, we consider the case $g=0$.
We argue our handle attachment.
\begin{itemize}
\item First we attach $k_1$ $2$-handles to the sphere disjointly and simultaneously. We must have a disjoint union of $k_1+1$ copies of the sphere $S^2$. 

More precisely, we choose $k_1$ disjoint copies of the product $S^1 \times D^1$ of a circle and a closed interval smoothly embedded in ${f_{0,v}}^{-1}(t_{v,1}+\epsilon)$, diffeomorphic to the sphere $S^2$.

We attach $2$-handles to the copies of $S^1 \times D^1$. The sphere ${f_{0,v}}^{-1}(t_{v,1}+\epsilon)$ is changed into a surface $S_{v,1}$, a disjoint union of $k_1+1$ copies of the sphere $S^2$.
\item After that we attach $k_2$ $1$-handles to change the previous surface $S_{v,1}$, consisting of $k_1+1$ copies of the sphere $S^2$ to a single $2$-dimensional sphere.
 
More precisely, we choose a $2$-dimensional disk $D^2$ smooothly embedded in the interior of each connected manifold which is a $2$-dimensional disk and represented as the intersection of each of the $k_1+1$ spheres from $S_{v,1}$ and the original sphere ${f_{0,v}}^{-1}(t_{v,1}+\epsilon)$. For $k_1-1$ of the $k_1+1$ connected components, we choose an additional $2$-dimensional disk $D^2$ smoothly embedded there and disjoint from the previously embedded disk.

We attach each $1$-handle to a disjoint union of two $2$-dimensional disks of these all $2k_2$ copies of $2$-dimensional disks to connect the $k_1+1$ copies of the sphere $S^2 \subset S_{v,1}$. We can connect the disjoint union $S_{v,1}$ of spheres to obtain a sphere as a result only in this way. In addition, the relation $k_1=k_2:=k$ must hold. 

\end{itemize}
Since the preimage ${q_{f_v}}^{-1}(v)$ must contain at least one singular point of it, for the STTF Morse function $f_{0,v}$ obtained by using a homotopy from $f_v$ and its handles, we also have $k_1=k_2=k>0$. 
We can also attach the handles disjointly and simultaneously. We can deform it by a suitable homotopy to a desired local Morse function $f$. Our desired Morse function must enjoy the property that at all singular points of it the values are same of course. The Reeb space $W_{f_{0,v}}$ is homeomorphic to a graph of 1st Betti number $k \geq 1$. 

We show an example for $k=2$ in FIGURE \ref{fig:2}.
\begin{figure}
		\includegraphics[width=80mm,height=40mm]{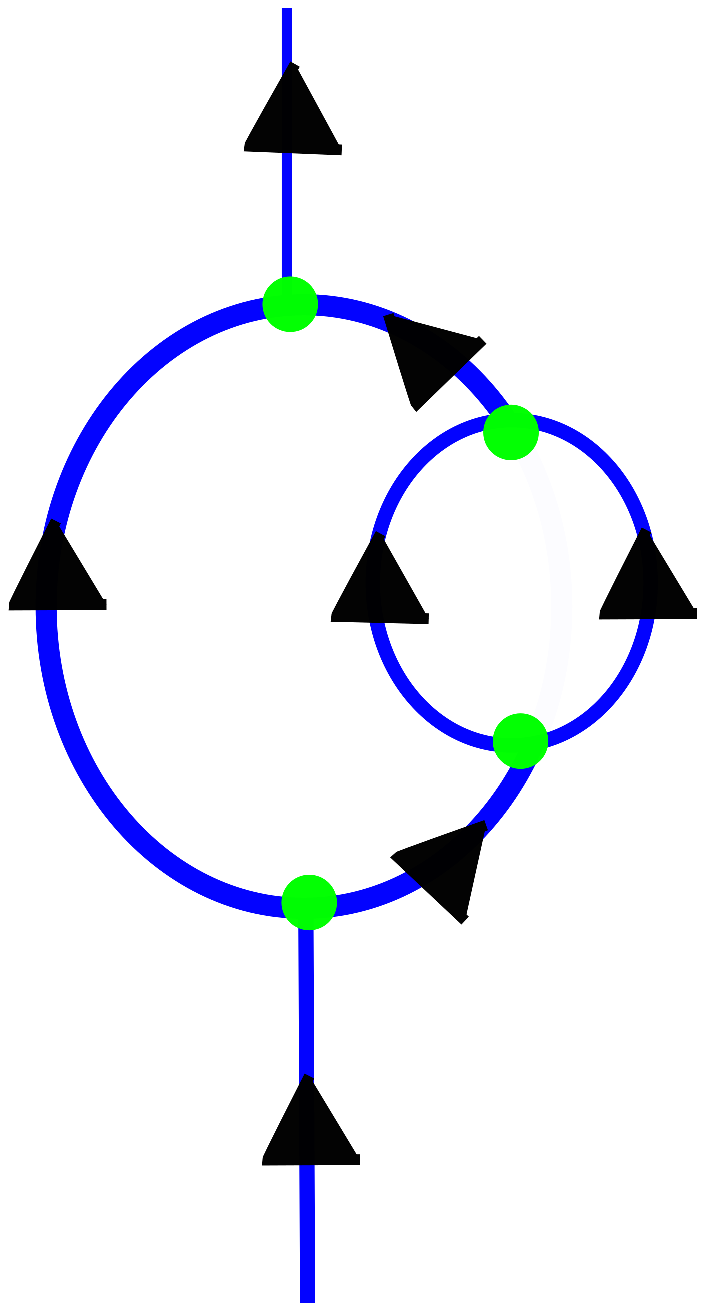}
		\caption{The Reeb digraph $W_{f_{0,v}}$ for $k=2$. The preimage of a point in the interior of a blue colored edge is diffeomorphic to $S^2$. Green dots are for singular points of index $1$ or $2$ for the Morse function.}
\label{fig:2}
\end{figure}
\item \label{3}
At each remaining vertex $v$, we can locally have a Morse function $f_v:X_v \rightarrow \mathbb{R}$ on a $3$-dimensional compact and orientable manifold $X_v$ enjoying the following properties as an example.
Let $i_{v,S^2} \geq 0$ ($i_{v,S^1 \times S^1} \geq 0$) denote the number of oriented edges entering $v$ where $l_{K}$ has the value $0$ (resp. $1$).   
Let $o_{v,S^2} \geq 0$ ($o_{v,S^1 \times S^1} \geq 0$) denote the number of oriented edges departing from $v$ where $l_{K}$ has the value $0$ (resp. $1$). 
The former case (\ref{2}) is for $i_{v,S^2}=o_{v,S^2}=1$ and $i_{v,S^1 \times S^1}=o_{v,S^1 \times S^1}=0$ and we do not consider this case here.
\begin{itemize}
\item At all singular points of the function the values are $g(v)$.
\item The image is a closed interval $[t_{v,1},t_{v,2}]$ containing $g(v)$ in its interior.
\item The preimage of $t_{v_1}$ is a disjoint union of $i_{v,S^2}$ copies of $S^2$ and $i_{v,S^1 \times S^1}$ copies of $S^1 \times S^1$.
\item The preimage of $t_{v_2}$ is a disjoint union of $o_{v,S^2}$ copies of $S^2$ and $o_{v,S^1 \times S^1}$ copies of $S^1 \times S^1$.
\item We can deform the Morse function by a suitable homotopy to an SSTF Morse function $f_{0,v}:X_v \rightarrow \mathbb{R}$ on our $3$-dimensional compact and connected manifold $X_v$ with the following properties.

\begin{itemize}
\item The number of singular points of the resulting SSTF Morse function $f_{0,v}$ is $2i_{v,S^1 \times S^1}+i_{v,S^2}+2o_{v,S^1 \times S^1}+o_{v,S^2}-2$. 
Let $t_{f_v,j}$ denote the $j$-th smallest value in the values realized as the values at singular points of the function. We also use the notation $t_{f_v,0}=t_{v,1}$ and $t_{f_v,2i_{v,S^1 \times S^1}+i_{v,S^2}+2o_{v,S^1 \times S^1}+o_{v,S^2}-2}=t_{v,2}$.
\item For each point $t \in (t_{f_v,j},t_{f_v,j+1})$ with $0 \leq j \leq i_{v,S^1 \times S^1}$, the preimage ${f_{0,v}}^{-1}(t)$ contains no singular point of $f_{0,v}$ and diffeomorphic to a disjoint union of $i_{v,S^2}+j$ copies of $S^2$ and $i_{v,S^1 \times S^1}-j$ copies of $S^1 \times S^1$.
\item For each point $t \in (t_{f_v,j},t_{f_v,j+1})$ with $i_{v,S^1 \times S^1}+1 \leq j \leq 2i_{v,S^1 \times S^1}+i_{v,S^2}-1$, the preimage ${f_{0,v}}^{-1}(t)$ contains no singular point of $f_{0,v}$ and diffeomorphic to a disjoint union of $i_{v,S^2}+2i_{v,S^1 \times S^1}-j$ copies of $S^2$.
\item For each point $t \in (t_{f_v,j},t_{f_v,j+1})$ with $2i_{v,S^1 \times S^1}+i_{v,S^2} \leq j \leq 2i_{v,S^1 \times S^1}+i_{v,S^2}+o_{v,S^2}+o_{v,S^1 \times S^1}-2$, the preimage ${f_{0,v}}^{-1}(t)$ contains no singular point of $f_{0,v}$ and diffeomorphic to a disjoint union of $j-2i_{v,S^1 \times S^1}-i_{v,S^2}+2$ copies of $S^2$.
\item For each point $t \in (t_{f_v,j},t_{f_v,j+1})$ with $2i_{v,S^1 \times S^1}+i_{v,S^2}+o_{v,S^2}+o_{v,S^1 \times S^1}-1 \leq j \leq 2i_{v,S^1 \times S^1}+i_{v,S^2}+2o_{v,S^1 \times S^1}+o_{v,S^2}-2$, the preimage ${f_{0,v}}^{-1}(t)$ contains no singular point of $f_{0,v}$ and diffeomorphic to a disjoint union of $2o_{v,S^2}+2o_{v,S^1 \times S^1}j+2i_{v,S^1 \times S^1}+i_{v,S^2}-j-2$ copies of $S^2$ and $j-2i_{v,S^1 \times S^1}-i_{v,S^2}-o_{v,S^2}-o_{v,S^1 \times S^1}+2$ copies of $S^1 \times S^1$.
\item The Reeb space $W_{f_{0,v}}$ is a graph of 1st Betti number $0$.
\end{itemize}
For a suitable Morse function explained just before, this can be also corresponded to attachment of handles to the surface
${f_{0,v}}^{-1}(t_{v,1}+\epsilon)$ to obtain $X_v={f_{0,v}}^{-1}([t_{v,1},t_{v,2}])$ according to the following procedure where $\epsilon>0$ is a sufficiently small number. 
\begin{itemize}
\item First we attach a $2$-handle to each of $i_{v,S^1 \times S^1}$ copies of the torus $S^1 \times S^1$ to change it to a copy of the sphere $S^2$. 
More precisely, for each torus, we choose a manifold $D_{v,1,j}$ diffeomorphic to $S^1 \times D^1$ which does not disconnect the torus in the torus. The manifold $D_{v,1,j}$ is indexed by an integer $1 \leq j \leq i_{v,S^1 \times S^1}$. 
We attach the handles to them.
We change the surface to another surface $S_{v,1}$, a disjoint union of $i_{v,S^2}+i_{v,S^1 \times S^1}$ copies of the sphere $S^2$.
\item Second we attach $i_{v,S^1 \times S^1}+i_{v,S^2}-1$ $1$-handles to change the previous surface to a copy of $S^2$. More precisely, for each of all $i_{v,S^2}+i_{v,S^1 \times S^1}$ disks represented as a connected component of the intersection of $S_{v,1}$ and ${f_{0,v}}^{-1}(t_{v,1}+\epsilon)$, we choose a smoothly embedded $2$-dimensional disk $D^2$ in its interior. Let $D_{v,2,j}$ denote the disk, indexed by an integer $1 \leq j \leq i_{v,S^2}+i_{v,S^1 \times S^1}$. For $i_{v,S^2}+i_{v,S^1 \times S^1}-2$ of all $i_{v,S^2}+i_{v,S^1 \times S^1}$ disks before, we add another smoothly embedded $2$-dimensional disk $D^2$ in its interior disjointly from the previously chosen disks $D_{v,2,j}$. Let $D_{v,2,j}$ denote the disk, indexed by an integer $i_{v,S^2}+i_{v,S^1 \times S^1}+1 \leq j \leq 2i_{v,S^2}+2i_{v,S^1 \times S^1}-2$. We attach the handles to them. We change the surface $S_{v,1}$ to another surface $S_{v,2}$, a copy of the $2$-dimensional sphere $S^2$.
\item Third we attach $o_{v,S^1 \times S^1}+o_{v,S^2}-1$ copies of $2$-handles to change the previous sphere $S_{v,2}$ to a disjoint union of $o_{v,S^1 \times S^1}+o_{v,S^2}$ copies of $S^2$. As before, we can also attach the handles disjointly and simultaneously on the interior of the intersection $S_{v,2} \bigcap {f_{0,v}}^{-1}(t_{v,1}+\epsilon)$. More precisely, we choose manifolds $D_{v,3,j}$ diffeomorphic to $S^1 \times D^1$ disjointly in such a way that they divide the sphere $S_{v,2}$ into $o_{v,S^1 \times S^1}+o_{v,S^2}$ copies of the $2$-dimensional sphere $S^2$. 
Let $S_{v,3}$ denote the resulting disjoint union of the $o_{v,S^1 \times S^1}+o_{v,S^2}$ copies of the $2$-dimensional sphere $S^2$.
The manifold $D_{v,3,j}$ is indexed by an integer $1 \leq j \leq o_{v,S^1 \times S^1}+o_{v,S^2}-1$. 
We also attach them disjointly from the $1$-handles and the $2$-handles of the previous two steps.
\item Last we choose $o_{v,S^1 \times S^1}$ spheres in the disjoint union $S_{v,3}$ of the previous spheres, attach a $1$-handle to each of the chosen sphere to make it diffeomorphic to the torus $S^1 \times S^1$. As before, we can also attach the handles disjointly and simultaneously on the interior of the intersection $S_{v,3} \bigcap {f_{0,v}}^{-1}(t_{v,1}+\epsilon)$, represented as a disjoint union of $o_{v,S^1 \times S^1}+o_{v,S^2}$ copies of the $2$-dimensional disk $D^2$. More precisely, we choose two disjoint copies $D_{v,4,j,j^{\prime}}$ ($j^{\prime}=1,2$) of the $2$-dimensional disk $D^2$ embedded smoothly in the interior of the $j$-th chosen disk ($1 \leq j \leq o_{v,S^1 \times S^1}$). We also attach them disjointly from the $1$-handles and the $2$-handles of the previous three steps.
\end{itemize}
\end{itemize}
We show an example for $i_{v,S^2}=i_{v,S^1 \times S^1}=o_{v,S^2}=o_{v,S^1 \times S^1}=1$ in FIGURE \ref{fig:3}. FIGURE \ref{fig:4} shows attachment of handles for this.
\begin{figure}
		\includegraphics[width=80mm,height=40mm]{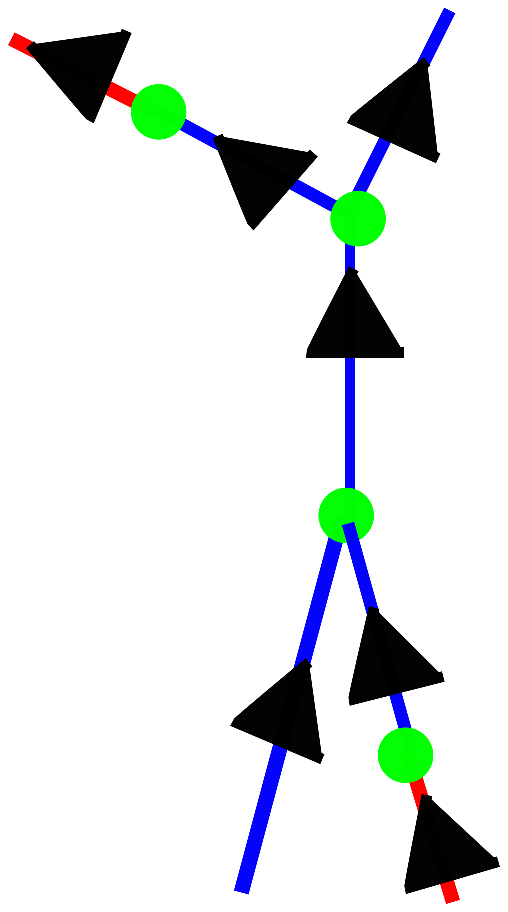}
		\caption{The Reeb digraph $W_{f_{0,v}}$ for $i_{v,S^2}=i_{v,S^1 \times S^1}=o_{v,S^2}=o_{v,S^1 \times S^1}=1$. The preimage of a point in the interior of a blue (red) colored edge is diffeomorphic to $S^2$ (resp. $S^1 \times S^1$). Green dots are for singular points of index $1$ or $2$ for the function.}
\label{fig:3}
	\end{figure}
\begin{figure}
		\includegraphics[width=80mm,height=40mm]{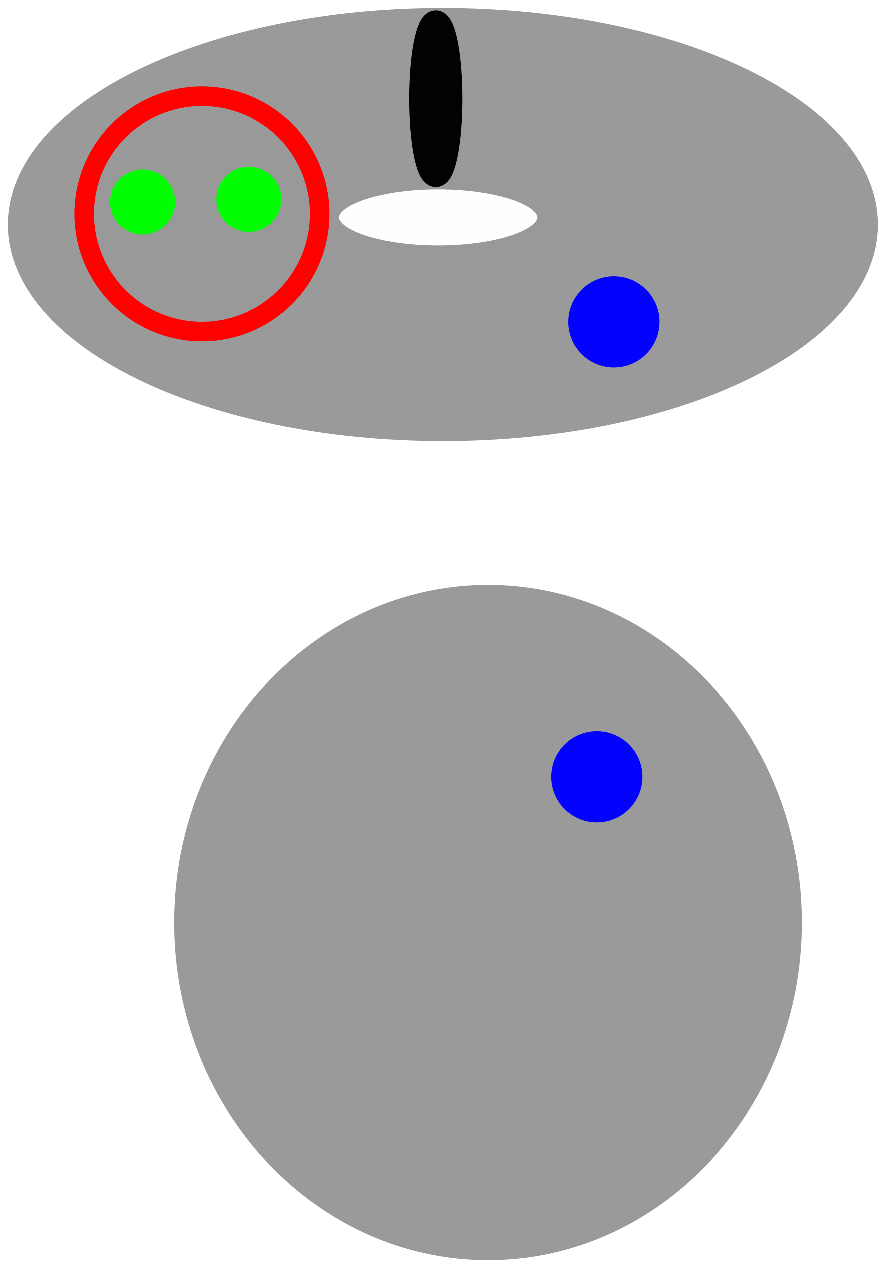}
		\caption{An example of attachment of handles to the surface ${f_{0,v}}^{-1}(t_{v,1}+\epsilon)$, colored in gray. This is also for the case $i_{v,S^2}=i_{v,S^1 \times S^1}=o_{v,S^2}=o_{v,S^1 \times S^1}=1$ or FIGURE \ref{fig:3}. The manifold $D_{v,1,1}$, diffeomorphic to $S^1 \times D^1$, is colored in black, the disks $D_{v,2,1}$ and $D_{v,2,2}$ are colored in blue, the manifold $D_{v,3,1}$, diffeomorphic to $S^1 \times D^1$, is colored in red, and the disks $D_{v,4,1,1}$ and $D_{v,4,1,2}$ are colored in green. To increase the 1st Betti number of a corresponding SSTF Morse function to $k>0$, we can add $k$ copies of the disk $D^2$ and $k$ copies of the product $S^1 \times D^1$ to which additional handles are attached to and they should be chosen apart from the existing copies of $D^2$ and $S^1 \times D^1$.}
\label{fig:4}
\end{figure}
\ \\
As the previous case, we can also increase $k \geq 0$ pairs of $1$-handles and $2$-handles and have a corresponding another Morse function suitably. We can also deform the Morse function by a suitable homotopy to an SSTF Morse function such that the Reeb space is homeomorphic to a graph of 1st Betti number $k \geq 0$. These additional handles must be and can be attached as in the previous case (\ref{2}). These handles can be also attached disjointly in submanifolds in the interior of a small disk apart from the handles, for example.  
\end{enumerate}

We give a short remark on this exposition on local functions by using handles. \cite{michalak1} respects this exposition using arguments on handles. \cite{kitazawa1} first respected this in earlier preprint versions. However, based on comments by referees and our considerations for example, we adopt another way, removing and attaching trivial smooth bundles whose fibers are diffeomorphic to $D^2$, in the published version. We do not explain this precisely. These two kinds of exposition are essentially same of course.

Around each edge, we construct a trivial smooth bundle over a closed interval in the interior of the edge. Its fiber is $S^2$ or $S^1 \times S^1$ according to the function $l_K$.
By gluing these local functions together, we have a desired Morse function $f:M \rightarrow \mathbb{R}$ on a $3$-dimensional closed, connected and orientable manifold $M$.
We can deform this into an SSTF Morse function $f_0:M \rightarrow \mathbb{R}$ by a suitable homotopy, which is also obtained by gluing each local SSTF Morse function $f_{0,v}$.   
We can see that the closed, connected and orientable manifold is diffeomorphic to a manifold obtained in the following way. We also present precise arguments.
\begin{itemize}
\item First we can see that for our resulting SSTF Morse function $f_0$ there exists exactly $c$ edges of $W_{f_0}$ such that the preimage ${q_{f_0}}^{-1}(t)$ of $t$ in the interior is diffeomorphic to the torus $S^1 \times S^1$. Such edges are not mutually adjacent by our construction. The 1st Betti number of $W_{f_0}$ is greater than or equal to $a+b$. We can make the 1st Betti number an arbitrary integer $r \geq a+b$ by our construction.
We remove the preimage of a small connected regular neighborhood of each of these $c$ edges. We have a $3$-dimensional manifold whose boundary is diffeomorphic to a disjoint union of $2c$ copies of $S^2$.
We attach copies of the cylinder $S^2 \times D^1$ instead to the boundary. This naturally yields another SSTF Morse function such that preimages of single points containing no singular points of it are disjoint unions of spheres. The resulting $3$-dimensional closed and connected manifold is diffeomorphic to the sphere $S^3$ in the case $a+b=0$ or a manifold diffeomorphic to a connected sum $\sharp (S^1 \times S^2)$ of $r \geq a+b$ copies of $S^1 \times S^2$. This is regarded as a higher dimensional version of Theorem \ref{thm:1}. We can understand the topology of the manifold $M$ from classical $3$-dimensional manifold theory for example. We may also have various ways of understanding. For example, on fundamental groups, we can also apply isomorphisms of fundamental groups of the manifold and the Reeb graph, discussed in \cite[Corollary 4.8]{saekisuzuoka}, \cite{kitazawa1}, and \cite[Corollary 4]{kitazawa2}.
\item The original manifold $M$ can be, thanks to fundamental arguments of $3$-dimensional manifolds, diffeomorphic to any manifold represented as a connected sum of $r \geq a+b$ copies of $S^1 \times S^2$ and $c^{\prime} \leq c$ manifolds whose fundamental groups are finite and whose Heegaard genera are $1$ (for any non-negative integer $c^{\prime} \leq c$). Furthermore, a $3$-dimensional closed, connected and orientable manifold admitting such an SSTF Morse function must be diffeomorphic to such a manifold, conversely. This is also important in STEP 3.
\end{itemize}

STEP 3 Topological restrictions on the $3$-dimensional closed and orientable manifold $M$ of the domain of the Morse function $f$. \\

Here we respect {\cite[Lemma 6.6]{saeki1}} and its original proof mainly: we introduce this in Theorem \ref{thm:5}.

We consider the case $g=1$ there. We also abuse most of the notation from STEP 2.

We can deform such a Morse function to an SSTF Morse function $f_0:M \rightarrow \mathbb{R}$ by a suitable homotopy. We can also have such a function $f_0$ in such a way that to each local SSTF Morse function $f_{0,v}: X_v \rightarrow \mathbb{R}$ around each vertex $v$, the following handle attachment can be corresponded. Let $f_{0,v}(X_v)=[t_{v,1},t_{v,2}]$ as in STEP 2.
\begin{enumerate}
\item First, to ${f_{0,v}}^{-1}(t_{v,1}+\epsilon)$, $2$-handles are attached. 
\item After that, to ${f_{0,v}}^{-1}(t_{v,1}+\epsilon)$, $1$-handles are attached. 
\end{enumerate}
This is from main ingredient of {\cite[Lemma 6.6]{saeki1}} and its original proof.  
We present new and further arguments.
 
From our construction of $f_0$ and the handle attachment, we have the following for the Reeb digraph $W_{f_0}$.
\begin{enumerate}
\item Consider the union $W_{f_0,S^1 \times S^1}$ of all edges of $W_{f_0}$ the preimage of each point in the interior of which is diffeomorphic to the torus $S^1 \times S^1$. FIGURE \ref{fig:1} with related arguments implies that the resulting set is (piecewisely smooth) homeomorphic to a disjoint union of closed intervals. 
\item For each connected component $I_v \subset W_{f_0,S^1 \times S^1}$ before, represented as a closed interval, the image $\bar{f_0}(I_v) \subset \mathbb{R}$ is not a subset of $[t_{v,1}+\epsilon,t_{v_2}-\epsilon]$.  
\end{enumerate}
This implies that the number of connected components of $W_{f_0,S^1 \times S^1}$ is at most $c$.

We can deform the function $f_0$ and the Reeb digraph $W_{f_0}$ one after another by local changes presented in FIGURE \ref{fig:5}, where colors are used as before. 

\begin{figure}
		\includegraphics[width=80mm,height=40mm]{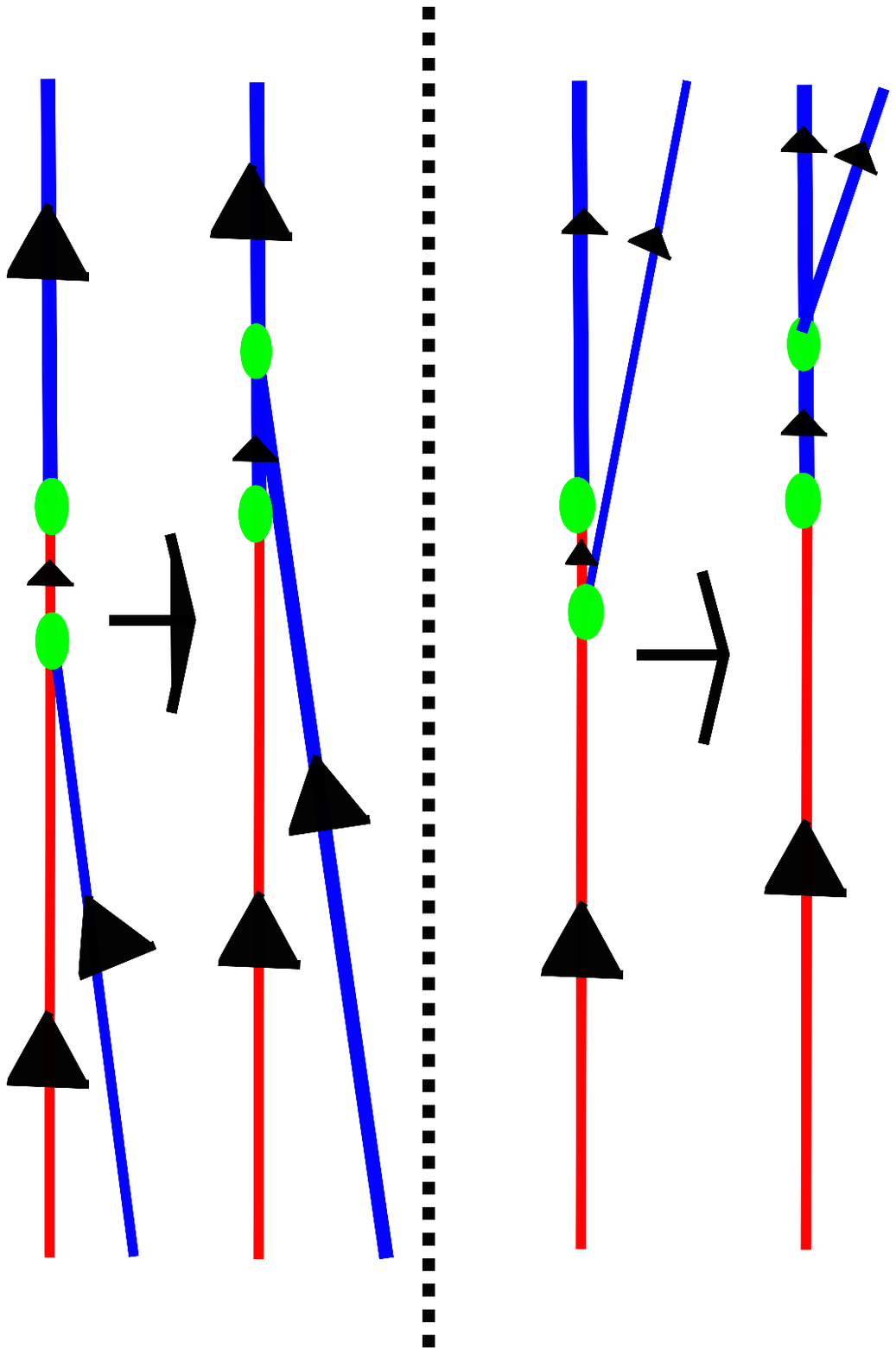}
		\caption{Important local changes of Reeb digraphs. Colors are used as in the previous Figures. Of course for these local SSTF functions $f_v$, we can consider $-f_v$ and have similar cases.}
\label{fig:5}
\end{figure}

This is thanks to the topology of the manifold and fundamental arguments on $3$-dimensional manifold theory. The $3$-dimensional manifolds of the preimages of the local graphs in FIGURE \ref{fig:5} are all diffeomorphic to $(S^1 \times D^2)-{\rm Int}\ (D^2 \sqcup D^2)$. Here ${\rm Int}\ X$ means the interior $X-\partial X$ of the topological manifold $X$, the two disjoint disks $D^2 \sqcup D^2$ are smoothly embedded in the interior of $S^1 \times D^2$ and here its interior is removed. Note that similar local changes are presented in \cite[Figure 4]{michalak2} and also used in \cite{marzantowiczmichalak} and that we cannot show this argument on the $3$-dimensional manifold from them directly. On the other hand, we can realize the local movements by applying suitable operations in \cite[Figure 5]{michalak2} twice. This gives another proof for the local movements. 

We have a new SSTF Morse function $f_0:M \rightarrow \mathbb{R}$ as discussed in the end of STEP 2. We can also know the type of the $3$-dimensional manifold $M$ as argued there. This completes STEP 3. \\
\ \\
This completes our proof.
\end{proof}
\section{Remarks.}
The following is an extension of Theorem \ref{thm:1}.

\begin{Thm}[\cite{gelbukh4}]
Let a connected digraph $K_g$ be defined from a piecewise smooth function $g:K \rightarrow \mathbb{R}$ with the condition that the degree of a vertex where the function $g$ has a local extremum is of degree at most $2$.
Then we have the following.
\begin{enumerate}
\item \label{thm:4.1} Suppose that at each vertex the function $g:K \rightarrow \mathbb{R}$ has a local extremum. Then we have a Morse-Bott function $f:M \rightarrow \mathbb{R}$ whose Reeb digraph $W_f$ and the digraph $K_g$ are isomorphic and $M$ must be diffeomorphic to the sphere $S^2$, the torus $S^1 \times S^1$, the projective plane {\rm (}a non-orientable surface of genus $1${\rm )}, or the Klein Bottle {\rm (}a non-orientable surface of genus $2${\rm )}. 

Conversely, for such a surface $M$ and such a digraph $K_g$, we can always have a Morse-Bott function $f:M \rightarrow \mathbb{R}$ whose Reeb digraph $W_f$ and the digraph $K_g$ are isomorphic.

\item \label{thm:4.2} Suppose that the previous condition {\rm (}\ref{thm:4.1}{\rm )} is not satisfied{\rm :} there exists a vertex where the function $g:K \rightarrow \mathbb{R}$ does not have a local extremum. 

Suppose that the 1st Betti number of the graph $K$ is $a$ and that the number of vertices of degree $2$ where $g_K$ does not have local extremums is $b$. 

Then for a closed, orientable and connected surface $M$ of genus $g \geq a+b$ or a closed, connected and non-orientable surface $M$ of genus $\max\{2a+b,1\}$, there exists a Morse-Bott function $f:M \rightarrow \mathbb{R}$ whose Reeb digraph $W_f$ and the digraph $K_g$ are isomorphic.

Conversely, for such a Reeb graph, if the domain of the Morse function is a closed and connected surface, then the surface is diffeomorphic to one of these surfaces. 

\end{enumerate}

\end{Thm}
 This is also extended further in \cite{gelbukh5} via several graph theoretic methods including theory founded in \cite{gelbukh3}.
\begin{Prob}
Can we extend our main result Theorem \ref{thm:2} to the Morse-Bott function case?
Related to this, can we extend Theorem \ref{thm:3} to the Morse-Bott function case?
\end{Prob}
We have a counterexample for Theorem \ref{thm:3}. Consider a connected digraph with exactly two edges connecting two distinct vertices: this is homeomorphic to the circle $S^1$.
We consider a bundle over the circle $S^1$ whose fiber is the torus $S^1 \times S^1$. We compose the projection of the bundle with the canonical projection of the circle to $\mathbb{R}$. This is a Morse-Bott
 function such that the preimages of points are diffeomorphic to the torus, or the disjoint union of two copies of the torus: at vertices, the preimages are diffeomorphic to the torus. The manifold of the domain is not diffeomorphic to one represented as a connected sum of finitely many copies of $S^1 \times S^2$ or manifolds whose fundamental groups are finite and whose Heegaard genera are $1$.

\begin{Prob}
Can we extend these studies to cases where genera of connected surfaces of preimages of single points containing no singular points may not be $0$ or $1$?
\end{Prob}
Related to this, we introduce \cite[Lemma 6.6]{saeki1}.
\begin{Thm}
\label{thm:5}
If a $3$-dimensional closed, connected and orientable manifold $M$
 admits a Morse function $f:M \rightarrow \mathbb{R}$ such that connected components of preimages of single points containing no singular points are closed surfaces whose genera are at most $g>0$, then $M$ admits such a Morse function $f_0:M \rightarrow \mathbb{R}$ with the additional property that at distinct singular points of $f_0$ the values are always distinct.

\end{Thm}
The cases $g=0,1$ are for the existence of an SSTF Morse function. This fact and related arguments are also important in the proof of Theorem \ref{thm:3} and our proof of Theorem \ref{thm:2}. However we do not know whether we can extend Theorem \ref{thm:3} to the desired cases.

Last, Michalak \cite{michalak3} independently studies Morse functions on $3$-dimensional closed, connected and orientable manifolds and their Reeb graphs where preimages of single points are not concerned. It mainly studies relations between the fundamental group of the manifold and the numbers of vertices of the Reeb graphs of degree $2$. As an explicit case, Theorem \ref{thm:6} has been obtained as its main result. We omit exposition on the Euler number of a circle bundle.
\begin{Thm}[{\cite[Theorem 1.1]{michalak3}}]
\label{thm:6}
For a circle bundle $M$ over a closed, connected and orientable surface of degree $g \geq 1$ with Euler number $\pm 1$, for any Morse function $f_{\rm m}:M \rightarrow \mathbb{R}$ with the following properties,
the 1st Betti number of the Reeb graph is $0$ and the number of vertices of degree $2$ is $4g$.
\begin{enumerate}
\item The function $f_{\rm m}:M \rightarrow \mathbb{R}$ is simple.
\item The number of singular points of the function $f_{\rm m}$ is equal to the minimum of the numbers of singular points of the functions in the set of all Morse functions on $M$ such that at distinct singular points of the functions the values are always distinct.
\end{enumerate}
\end{Thm} 
Furthermore, Michalak also studies a kind of characterization of a connected sum of the form $\sharp (S^2 \times S^1)$ and its connected sum with a single Len space by SSTF Morse functions in our present study, independently. Note that Michalak does not introduce, name, or use the notion of SSTF Morse function there. See \cite[Proposition 4.2 and Theorem 4.3]{michalak3}.
\section{Acknowledgement.}
The author would like to thank members of the research group of Osamu Saeki. The author would also like to thank people organizing and supporting Saga Souhatsu Mathematical Seminar (http://inasa.ms.saga-u.ac.jp/Japanese/saga-souhatsu.html): the author would also like to thank Inasa Nakamura for inviting the author as a speaker and letting him present \cite{kitazawa3, kitazawa5}. Discussing several Morse functions and topological properties of the manifolds of their domains with these persons and groups has motivated the author to challenge the present study. 
\section{Conflict of interest and Data availability.}
\noindent {\bf Conflict of interest.} \\
The author works at Institute of Mathematics for Industry (https://www.jgmi.kyushu-u.ac.jp/en/about/young-mentors/) and this study is closely related to our study. Our study thanks them for the supports. The author is also a researcher at Osaka Central
Advanced Mathematical Institute (OCAMI researcher), supported by MEXT Promotion of Distinctive Joint Research Center Program JPMXP0723833165. Although he is not employed there, our study also thanks them. \\
Some of works by other researchers and this version may overlap in some of the contents due to the nature that our problems are natural in theory of Morse functions and applications to differential topology and that related mathematical studies are very fundamental and classical in some senses, for example. However the present version of our paper is presented independent of these work. \\
Saga Souhatsu Mathematical Seminar (http://inasa.ms.saga-u.ac.jp/Japanese/saga-souhatsu.html), inviting the author as a speaker, is supported by JST Fusion Oriented REsearch for disruptive Science and Technology JPMJFR202U and the author was a speaker on 2024/7/12 supported by this.\\
\ \\
{\bf Data availability.} \\
Data essentially supporting our present study are all in the paper.

\end{document}